\documentclass[12pt,a4paper]{article}
\usepackage[latin1]{inputenc}
\usepackage[T1]{fontenc}
\usepackage[english, french]{babel}
\usepackage{amssymb,amsmath}
\usepackage[french]{minitoc}
\usepackage{epsfig,graphics}
\begin{document}
\newenvironment{dem}{\textbf{Proof}}{}
\newenvironment{rem}{\textit{\textbf{Remark:}}}{}
\newenvironment{remf}{\textit{\textbf{Remarque:}}}{}
\newcommand{\dmath}{\displaystyle}
\newtheorem{theo}{Theorem}[section]
\newtheorem{theof}{Théorème}[section]
%[section]
\newtheorem{prop}{Proposition}[section]
%[section]
\newtheorem{coro}{Corollary}[section]
\newtheorem{corof}{Corollaire}[section]
%[section]
\newtheorem{lem}{Lemma}[section]
\newtheorem{lemf}{Lemme}[section]
%[section]
\newtheorem{conj}{Conjecture}[section]
\newenvironment{abst}{\textbf{abstract}}{}
\author{ Jean-Christophe Bourgoin }
%
%
%\newline
%\newline
%
\date{}
%\newline
%
\title{The minimality of the map $\frac{x}{\|x\|}$ for weighted
energy.} \maketitle
\textsc{\small{Laboratoire de Mathematiques et Physiques
Théorique, Université de Tours, Parc Grandmont 37200 TOURS}}
\newline
\newline
\begin{abst}

In this paper, we investigate the minimality of the map
$\frac{x}{\|x\|}$ from the euclidean unit ball $\mathbf{B}^n$ to
its boundary $\mathbb{S}^{n-1}$ for weighted energy functionals of
the type $E_{p,f}= \int_{\mathbf{B}^n}f(r)\|\nabla u\|^p dx$,
where $f$ is a non-negative function. We prove that in each of the
two following cases:\\ i) $p=1$ and $f$ is non-decreasing,
\\ i)) $p$ is  an integer, $p \leq n-1$ and $f= r^{\alpha}$ with
$\alpha \geq 0$,\\ the map  $\frac{x}{\|x\|}$ minimizes $E_{p,f}$
among the maps in  $W^{1,p}(\mathbf{B}^n, \mathbb{S}^{n-1})$ which
coincide with  $\frac{x}{\|x\|}$  on $\partial \mathbf{B}^n$. We
also study the case where $ f(r)= r^{\alpha}$ with $-n+2 < \alpha
< 0$ and prove that   $\frac{x}{\|x\|}$ does not minimize
$E_{p,f}$ for $\alpha$ close to $-n+2$ and when $n \geq 6$, for
$\alpha$ close to $4-n$.
\end{abst}
\newline
\newline
\textbf{\small{Keys Words}}: minimizing map, $p$-harmonic map,
$p$-energy, weighted energy.
%
%
%
%\addcontentsline{toc}{subsection}{1.1.1 Introduction and statement
%of results .}
%
%
\subsection{ Introduction and statement of results}
For $n \geq 3$, the map $u_0(x)\!=\!\frac{x}{\|x\|}: \mathbf{B}^n
\longrightarrow \mathbb{S}^{n-1}$ from the unit ball
$\mathbf{B}^n$ of $ \mathbb{R}^n$ to its boundary
$\mathbb{S}^{n-1}$ plays a crucial role in the study of certain
natural energy functionals. In particular, since the works of
Hildebrandt, Kaul and Widman (\cite{HKW}), this map is considered
as a natural candidate to realize, for each real number $p \in
[1,n)$ the minimum of the $p$-energy functional,
\[
E_p(u)= \int_{\mathbf{B}^n} \|\nabla u\|^p dx
\]
among the maps $u \in
W^{1,p}(\mathbf{B}^n,\mathbb{S}^{n-1})\!=\!\{u \in
W^{1,p}(\mathbf{B}^n,\mathbb{R}^n; \|u\|\!=\!1\, a.e.\}$
satisfying $u(x)=x$ on $\mathbb{S}^{n-1}$.\\
This question was first treated in the case $p=2$. Indeed, the
minimality of $u_0$ for $E_2$ was etablished by Jäger and Kaul (
\cite{JK}) in dimension $n \geq 7$ and by Brezis, Coron and Lieb
in dimension 3 ( \cite{BCL}). In \cite{CG}, Coron and Gulliver
proved the minimality of $u_0$ for $E_p$ for any integer $p\in
\{1,\cdots,n-1\}$ and any dimension $n \geq 3$.\\

Lin (\cite{L}) has introduced the use of the elegant null
Lagrangian method (or calibration method) in this topic.
Avellaneda and Lin showed  the efficiency of this method in
\cite{AL} where they give a simpler alternative proof to the
Coron-Gulliver result. Note that several results concerning the
minimizing properties of $p$-harmonic diffeomorphisms  were also
obtained in this way in particular by Coron, Helein and El Soufi,
Sandier ( \cite{CHe}, \cite{He}, \cite{ES} and \cite{EJ}).

The case of non-integer $p$ seemed to be rather difficult. It is
only ten years after the Coron-Gulliver article \cite{CG}, that
Hardt, Lin and Wang (\cite{HLW1}) succeeded to prove that, for all
$n \geq 3$, the map $u_0$ minimizes $E_p$ for  $p \in [n-1,n)$.
Their proof is based on a deep studies of singularities of
harmonic and minimizing maps made in the last two decades. In
dimension $n \geq 7$, Wang (\cite{W}) and Hong (\cite{Ho1}) have
independently proved the minimality of $u_0$ for any $p\geq 2$
satisfying $p +2 \sqrt{p} \leq n-2$.

In \cite{Ho2}, Hong remarked that the minimality of the $p$-energy
$E_p$, $p \in (2, \nolinebreak[4]n-\nolinebreak[4]1]$, is related
to the minimization of the following weighted 2-energy:
\[
\tilde{E}_p(u)= \int_{\mathbf{B}^n} r^{2-p}\|\nabla u\|^2 dx
\]
where $r=\|x\|$. Indeed, using Hölder inequality, it is easy to
see that if the map $u_0$ minimizes $\tilde{E}_p$, then it also
minimizes  $E_p$ (see \cite{Ho2}, p.465). Unfortunately, as we
will see in Corollary 1.1 below, for many values of $p \in (2,n)$,
the map $u_0$ is not a minimizer of $\tilde{E}_p$. Therefore,
Theorem 6 of (\cite{Ho2}), asserting that $u_0$ minimizes
$\tilde{E}_p$ seems to be not correct and the question of whether
$u_0$ is a minimizing map of the $p$-energy $E_p$ for non-integer
$p \in (2, n-1)$ is still open \footnote{We suspect a problem in
Theorem 6 p.464 of \cite{Ho2}. Indeed the author claims that the
quantity $G_{\varphi_1^0, \cdots, \varphi_{n-1}^0}(v,p)$, which
represents a weighted energy of the map $v$ on the $3$-dimensional
cone $\mathcal{C}_0$ in $\mathbf{B}^n$, is uniformly proportional
to the weighted energy on the euclidian ball $\mathbf{B}^3$. There
is no reason for this fact to be true, the orthogonal projection
of $\mathcal{C}_0$ on to $\mathbf{B}^n$ being not homothetic.}

The aim of this paper is to study the minimizing properties of the
map $u_0$ in regard to some weighted energy functionals of the
form:
\[
E_{p,f}(u)= \int_{\mathbf{B}^n} f(r)\|\nabla u\|^p dx,
\]
where $p \in \{1, \cdots, n-1\}$ and $f\!:\! [0,1] \rightarrow
\mathbb{R}$ is a non-negative non-decreasing continuous function.
For $p=1$, the map $u_0$ minimizes $E_{1,f}$ for a large class of
weights. Indeed, we have the following
\begin{theo}
Suppose that $f$ is a non-negative differentiable non-decreasing
function. Then the map $u_0=\frac{x}{\|x\|}$ is a minimizer of the
energy $E_{1,f}$, that is, for any $u$ in
$W^{1,1}(\mathbf{B^n},\mathbb{S}^{n-1})$ with $u(x)=x$ on
$\mathbb{S}^{n-1}$, we have
\[
\int_ {\mathbf{B}^n} f(r) \|\nabla u_0\| dx \leq \int_
{\mathbf{B}^n} f(r) \|\nabla u\| dx,
\]
Moreover, if $f$ has  no critical points in $(0,1)$, then the map
$u_0=\frac{x}{\|x\|}$ is the unique minimizer of the energy
$E_{1,f}$, that is, the equality in the last inequality holds if
and only if $u=u_0$.\\
\end{theo}
For $p \geq 2 $, we restrict ourselves to power functions
$f(r)=r^{\alpha}$,
\begin{theo}
For any $\alpha \geq 0$ and any integer $p \in \{1, \cdots,
n-1\}$, the map $u_0=\frac{x}{\|x\|}$ is a minimizer of the energy
$E_{p,r^{\alpha}}$ that is, for any $u$ in
$W^{1,p}(\mathbf{B^n},\mathbb{S}^{n-1})$ with $u(x)=x$ on
$\mathbb{S}^{n-1}$, we have,
\[
\int_ {\mathbf{B}^n} r^{\alpha} \|\nabla u_0\|^p dx \leq \int_
{\mathbf{B}^n} r^{\alpha} \|\nabla u\|^p dx \,.
\]
Moreover, if $\alpha >0$, then the map $u_0= \frac{x}{\|x\|}$ is
the unique minimizer of the energy $E_{p,r^{\alpha}}$, that is the
equality in the last inequality holds if and only if $u=u_0$.
\end{theo}
The proof of these two theorems is given in section 2. It is based
on a construction of an adapted null-Lagrangian. The case of $p=1$
can be obtained passing through  more direct ways and will be
treated independently.\\

The case of weights of the form $f(r)= r^{\alpha}$, with
$\alpha<0$, is treated in section 3. The weighted energy
$\int_{\mathbf{B}^n}r^{\alpha}\|\nabla u_0\|^2 dx$ of
$u_0=\frac{x}{\|x\|}$ is finite for $\alpha > -n+2$. Hence we
consider the family of maps,
\[
u_a(x)= a+ \lambda_a(x)(x-a), \quad a \in \mathbf{B}^n,
\]
where $\lambda_a(x)\!\in\!\mathbb{R}$ is chosen such that
$u_a(x)\in \mathbb{S}^{n-1}$ (that is $u_a(x)$ is the intersection
point of $\mathbb{S}^{n-1}$ with the half-line of origin $a$
passing by $x$).
%
%
%
%\newpage
%
%\vspace{1cm}

%\newpage
%
We study the energy $E_{2,r^{\alpha}}(u_a)$ of these maps and
deduce the following theorem.
\begin{theo}
Suppose that $n\geq 3$.
\newline (i)\, For any $a \in \mathbf{B}^n
, a \neq 0$ , there exists a negative real number\\ $\alpha_0 \in
(-n+2,0)$, such that, for any  $\alpha \in (\!-\!n+2,\alpha_0]$ we
have
\[
\int_ {\mathbf{B}^n} r^{\alpha} \|\nabla u_0\|^2dx > \int_
{\mathbf{B}^n} r^{\alpha} \|\nabla u_a\|^2 dx \quad.
\]
(ii)\,For any integer $n \geq 6$, there exists $\alpha_0 \in
(4-n,5-n)$ such that, for any $\alpha \in (4-n, \alpha_0)$, there
exists $a \in \mathbf{B}^n$ such that,
\[
\int_ {\mathbf{B}^n} r^{\alpha} \|\nabla u_0\|^2dx > \int_
{\mathbf{B}^n} r^{\alpha} \|\nabla u_a\|^2 dx \quad.
\]
\end{theo}
Replacing in Theorem 0.3 $\alpha$ by $2-p$, $p \in (2,n)$, we
obtain the following corollary:
\begin{coro}
For any $n \geq 6$, there exists $p_0 \in (n-3,n-2)$ such that,
for any $p \in (p_0, n-2)$ the map $u_0= \frac{x}{\|x\|}$ does not
minimize the functional $\int_ {\mathbf{B}^n} r^{2-p} \|\nabla
u\|^2dx$ among the maps $u \in
W^{1,2}(\mathbf{B}^n,\mathbb{S}^{n-1})$ satisfying $u(x)=x$ on
$\mathbb{S}^{n-1}$.
\end{coro}
\textit{acknowledgements}. The author would express his grateful
to Professor Ahmad El Soufi and Professor Etienne Sandier for
their helpful advice.

%
%
%
%\addcontentsline{toc}{subsection}{1.1.2 proof of theorems 1.1 and
%1.2.}
%

\subsection{ Proof of theorems 0.1 and 0.2}
Consider an integer $p \in \{1,\cdots\!,n\!-\!1\}$ \, and $f$ a
differentiable, non-negative, increasing, and non-identically zero
map. We can suppose without loss of generality, that
$f(1)=1$.\\For any subset $I=\{i_1,\cdots\!, i_p \} \subset
\{1,\cdots\!,n\!-\!1\}$\, with $i_1\!<\!i_2\ldots<\!i_p$  and for
any map,
\[
u=(u_1,\cdots\!,u_n):\mathbf{B}^n \longrightarrow \mathbb{S}^{n-1}
\quad \textrm{in} \,\, \mathcal{C}^{\infty}(\mathbf{B}^n
,\mathbb{S}^{n-1}) \quad \textrm{with}\,\ u(x)=x\,\, \text{on}\,\,
\mathbb{S}^{n-1},
\]
we consider the n-form:
\[
\omega_I (u)= dx_1\wedge \cdots \wedge d(f(r) u_{i_1})\wedge\cdots
\wedge d(f(r) u_{i_k}) \wedge\cdots
 \wedge dx_n
\]
\begin{lem}
We have the identity:
\[
\int_{\mathbf{B}^n} \omega_I (u)=\,\int_{\mathbf{B}^n} \omega_I
(Id) \quad \forall \, x \in \mathbf{B}^n \,\quad \textrm{where}
\quad Id(x)=x.
\]
\end{lem}
\begin{dem}
By Stokes theorem, we have:
\begin{eqnarray*}
\int_{\mathbf{B}^n} \omega_I (u)
                             &=&  \int_{\mathbf{B}^n}  dx_1\wedge \cdots
\wedge d(f(r) u_{i_1})\wedge\cdots \wedge d(f(r) u_{i_p})
\wedge\cdots \wedge dx_n\\
 &=& \int_{\mathbf{B}^n} (-1) ^{i_1\!-\!1} d
\Big(\!f(r) u_{i_1} dx_1\wedge \cdots \wedge \widehat{d(f(r)
 u_{i_1})} \wedge\\
&& \hspace{5cm}\cdots \wedge\! d(f(r) u_{i_p})\! \wedge\cdots
\wedge dx_n \Big)
\\
 &=& \int_{\mathbb{S}^{n-1}} (-1)
^{i_1-1} x_{i_1} dx_1\wedge \cdots \wedge \widehat{d(f(r)
 u_{i_1})} \wedge\\
 && \hspace{5cm}\cdots \wedge d(f(r) u_{i_p})
\wedge\cdots \wedge dx_n .\\
\end{eqnarray*}
Indeed, on $\mathbb{S}^{n-1}$, we have  $ f(r)u_{i_1}=x_{i_1}$
($r=1, f(1)=1 \,\,\text{and}\,\,u(x)=x$). Iterating, we get the
designed identities.
\end{dem}
Consider the n-form:
\begin{displaymath}
S(u)= \sum_{|I|=p} w_I(u)
\end{displaymath}
By Lemma 0.1, we have:
\[
\int_{\mathbf{B}^n} S(u) =
\displaystyle\sum_{|I|=p}\int_{\mathbf{B}^n} w_I(u)= \sum_{|I|=p}
\int_{\mathbf{B}^n} d x=  C^{p}_{n} \frac{|\mathbb{S}^{n-1}|}{n},
\]
where $|\mathbb{S}^{n-1}|$ is the Lebesgue measure of the sphere.
\begin{lem}
The n-form $S(u)$ is $O(n)\!-\! equivariant$, that is, for any
rotation $R$ in $O(n)$, we have:
\[
S(^t R  u R) (^t Rx)= S(u) (x) \quad \forall x \in \mathbf{B}^n.
\]
\end{lem}
\begin{dem}
Consider  $S(u)(x)(e_1,\ldots, e_n)$ where $(e_1, \ldots, e_n)$ is
the stantard basis of $\mathbb{R}^n$ and notice that it is  equal
to $(-1)^n$ times  the $(p+1)^{th}$ coefficient of the polynomial
$P(\lambda)= \det(Jac(fu)(x)- \lambda Id)$ which does not change
when we replace $fu$ by $^tRfuR$.
\\
\end{dem}\\
For any $x \in \mathbf{B}^n$, let  $R\!\in\!  O(n)$ be such that $
^t Ru(x)\! =\! e_n=\! (0,\ldots,0,1)\!$. Consider  $y= ^t\!R
x\,\,,  v=^t\! RuR$, so that:
\[
v(y) =e_n, \quad  d(^t RuR )(y) (\mathbb{R}^n) \subset
e_n^{\bot}\quad \textrm{that is}\quad \frac{\partial v_n}{\partial
x_j}(y)=0 \quad \forall j\in\{1,\cdots,n\}.
\]
\begin {lem}
Let $a_1,\ldots\!,a_n$ be n non-negative numbers, and $p \in \{1,
\ldots,n- \nolinebreak[4] 1\}$. Then:
\[
\dmath\sum _{i1<\ldots<i_p}  a_{i_1} \cdots a_{i_p} \leq
\frac{1}{(n-1)^p} C_{n-1}^p \Big(\dmath\sum_{j=1}^{n-1} a_j
\Big)^p.
\]
\end{lem}
\begin{dem}
See for instance Hardy coll.[4], theorem 52.
\\
\end{dem}\\
Let $I= \{i_1,\cdots\!,i_p\}  \subset \{1,\cdots\!,n \}$. We
have:\\ if $i_p \neq n$,
\begin{eqnarray*}
\omega_I(v)(y)&=& \big(
 dx_1\wedge \cdots \wedge d(f(r) v_{i_1})\wedge\cdots
\wedge d(f(r) v_{i_k}) \wedge\cdots \wedge dx_n \big)(y)\\ &=&
|f(r)|^p \big(dx_1\wedge \cdots \wedge d v_{i_1}\wedge\cdots
\wedge d v_{i_k} \wedge\cdots \wedge dx_n \big)(y).
\end{eqnarray*}
Indeed, $\forall j \leq n\!-\!1$, $d(f(r) v_j(y))= d(f(r)) v_j(y)+
f(r) dv_j(y)=\!f(r)dv_j(y)$ since $v(y)\!=\!e_n$.\\
If $i_p=n$,
\begin{eqnarray*}
\omega_I(v)(y)&=& |f(r)|^{p-1} \big(
 dx_1\wedge \cdots \wedge d v_{i_1}\wedge\cdots
\wedge df  \big)(y).
\end{eqnarray*}
Indeed, $d(f(r)v_n)(y)= df(y) v_n (y)+ f(r)dv_n(y)= df(y)$ (as
$dv(y) \subset \nolinebreak[4] e_n^{\perp}$).
The Hadamard inequality gives:
\begin{eqnarray*}
|S(v)(y)|=
\big|\dmath\sum_{|I|=p} \omega_I(v)(y)\big|
&\leq&
|f(r)|^p \dmath\sum_{1\leq i_1 <i_2<\ldots <i_p\leq n-1}\!\!
\|dx_1\|\cdots \|d v_{i_1}\|\\
&& \hspace{4cm}\cdots \|dv_{i_p}\| \cdots \|dx_n\| (y)\\
&+&
|f(r)|^{p-1} \!\!\dmath\sum_{1\leq i_1 <i_2<\ldots <i_{p-1}\leq
n-1} \!\!\|dx_1\|\cdots \|dv_{i_1}\|\\
&&
\hspace{4cm} \cdots \|dv_{i_p}\| \cdots \!\|df\| \!(y)\\
&\leq &
|f(r)|^p  \bigg(\!\dmath\sum_{1\leq i_1 <i_2<\ldots <i_p\leq
n-1}\!\!\! \|dx_1\|^2\cdots \|d v_{i_1}\|^2\!\cdots\\
&&
\quad \quad \quad \cdots \|dv_{i_p}\|^2 \!\cdots \|dx_n\|^2
(y)\!\bigg)^{\frac{1}{2}}\big(C_n^p\big)^{\frac{1}{2}}\\
&+ &
f'(r)f(r)^{p-1}\!\!\dmath\sum_{1\leq i_1 <i_2<\ldots <i_{p-1}\leq
n-1}\!\! \|dx_1\|\cdots \|dv_{i_1}\|\\
&&
\hspace{4cm}\cdots \|dv_{i_p}\|(y).
\end{eqnarray*}
The Hardy inequality gives, after integration and using the fact
that $\|\nabla u\|=\|\nabla v\|$,
\begin{eqnarray*}
\frac{C_n^p}{n}|\mathbb{S}^{n-1}| &\leq&
\frac{C_{n-1}^p}{(n-1)^{p/2}} \int_{\mathbf{B}^n} f^p(r)
 \|\nabla
u\|^p dx\\
&&
+ \frac{C_{n-1}^{p-1}}{(n-1)^{\frac{p-1}{2}}} \int_{\mathbf{B}^n}
f'(r)f^{p-1}(r) \|\nabla u\|^{p-1}dx.\quad (1)
\end{eqnarray*}
\begin{rem}:\,
If $f'$ is positive and if equality  holds in (1), then, $\forall
i \leq n-1$, $y_i=0$ and $y_n= \pm \frac{x}{\|x\|}$, which implies
that $u(x)= \pm \frac{x}{\|x\|}$.
%
%The case of equality in inequality of Hadamard
%gives that if $\forall i \leq n-1,\, \frac{\partial f}{\partial
%x_i}(y)=0$ therefore $ \forall i \leq n-1,\,y_i=0$. As
%$\|y\|=\|x\|$, we deduce that $y_n=\|x\|$ or $y_n=-\|x\|$ then, if $f'$ does not vanish, the inequality (1) becomes a equality
%if and only if  $u(x)=\pm \frac{x}{\|x\|}$
\end{rem}\\
\newline
\begin{dem}
\textbf{of the Theorem 1.1}
Inequality (1) give
\[
|\mathbb{S}^{n-1}| \leq  \sqrt{n-1} \int_{\mathbf{B}^n} f(r)
 \|\nabla u\|dx + \int_{\mathbf{B}^n} f'(r) dx.
\]
Hence:
\[
 \int_{\mathbf{B}^n} f
 \|\nabla u\|dx\geq
\frac{|\mathbb{S}^{n-1}|}{\sqrt{n-1}} \big(1- \int_0^1
f'(r)r^{n-1} dr \big)
\]
\[
 \int_{\mathbf{B}^n} f
 \|\nabla u\|dx\geq
\sqrt{n-1}|\mathbb{S}^{n-1}|  \int_0^1 f(r)r^{n-2}
dr=\int_{\mathbf{B}^n} f(r)
 \|\nabla
u_0\|dx.
\]
To see the uniqueness il suffices to refer to the remark above. It
gives that for any $x \in \mathbf{B}^n$, $u(x)=\frac{x}{\|x\|}$ or
$u(x)=-\frac{x}{\|x\|}$. As $u(x)=x$ on the unit sphere, we have,
for any $x \in \mathbf{B}^n \backslash \{0\}$,
$u(x)=\frac{x}{\|x\|}$. \quad \quad \quad \quad  $\blacksquare$\\
\end{dem}
\newline
\begin{dem}
%\label{2.4}
\textbf{of the Theorem 1.2}. Let $\alpha$ be a positive real
number. From inequality (1) we have:
\[
\frac{C_n^p}{n}|\mathbb{S}^{n-1}| \leq
\frac{C_{n-1}^p}{(n-1)^{p/2}} \int_{\mathbf{B}^n} r^{\alpha p}
\|\nabla u\|^p dx +\alpha
\frac{C_{n-1}^{p-1}}{(n-1)^{\frac{p-1}{2}}} \int_{\mathbf{B}^n}
r^{\alpha p-1} \|\nabla u\|^{p-1} dx.
\]
By Hölder inequality, we have, setting $q=\frac{p}{p-1}$:
\begin{eqnarray*}
\frac{C_n^p}{n}|\mathbb{S}^{n-1}| &\leq&
\frac{C_{n-1}^p}{(n-1)^{p/2}} \int_{\mathbf{B}^n} r^{\alpha p}
 \|\nabla
u\|^p dx \\ && +\quad \alpha
\frac{C_{n-1}^{p-1}}{(n-1)^{\frac{p-1}{2}}} \!\left(
\int_{\mathbf{B}^n}\!r^{ p(\alpha-1)}dx
\!\right)^{1/p}\!\!\!\left(\! \int_{\mathbf{B}^n} \!\!r^{\alpha p}
\|\nabla u\|^p dx\!\right)^{1/q}\\ &\leq&
\frac{C_{n-1}^p}{(n-1)^{p/2}} \int_{\mathbf{B}^n} r^{\alpha p}
 \|\nabla
u\|^p dx\\ && +\quad \alpha
\frac{C_{n-1}^{p-1}}{(n-1)^{\frac{p-1}{2}}}\frac{|\mathbb{S}^{n-1}|^{1/p}}{(n+
p(\alpha-1))^{1/p}} \left(\int_{\mathbf{B}^n} \!r^{\alpha p}
 \|\nabla
u\|^p dx \right)^{1/q}
\end{eqnarray*}
Consider  the polynomial function:
\[
P(t)=\frac{C_{n-1}^p}{(n-1)^{p/2}} t^q + \alpha
\frac{C_{n-1}^{p-1}}{(n-1)^{\frac{p-1}{2}}}\frac{|\mathbb{S}^{n-1}|^{1/p}}{(n+
p(\alpha-1))^{1/p}}t-\frac{C_n^p}{n}|\mathbb{S}^{n-1}|.
\]
Setting $A=\left(\int_{\mathbf{B}^n} r^{\alpha p}
 \|\nabla
u\|^p \right)^{1/q}$\,\,and $B=\left(\int_{\mathbf{B}^n} r^{\alpha
p}
 \|\nabla
u_0\|^p \right)^{1/q}$, we get $ P(A)\geq 0$ while
\begin{eqnarray*}
P(B) &= &\frac{C_{n-1}^{p-1}}{n+
p(\alpha-1)}|\mathbb{S}^{n-1}|+\alpha \frac{C_{n-1}^{p-1}}{n+
p(\alpha-1)}|\mathbb{S}^{n-1}|
-\frac{C_n^p}{n}|\mathbb{S}^{n-1}|\\ &=& \frac{C_{n-1}^{p-1}}{n+
p(\alpha-1)}|\mathbb{S}^{n-1}|\left(\frac{n-p}{n}+\alpha
-\frac{C_n^p}{nC_{n-1}^{p-1}}(n+ p(\alpha-1))\right)\\ &=& 0.
\end{eqnarray*}
On the other  hand, $ \forall t \geq 0,\,  P'(t)>0$. Hence, $P$ is
increasing in $[0,+ \infty)$ and is equal to zero only for $B$.
Necessarily, we have $A\!\geq \nolinebreak[4] B$.\\ Moreover, if
$\alpha > 0$, $A=B$ implies that equality in the inequality (1)
holds. Referring to the remark above, and as $u_0(x)=x$ on the
sphere, we have $u=u_0= \frac{x}{\|x\|}$. Replacing $ \alpha$ by $
\alpha / p $ we finish the prove of the theorem. \quad\quad
\quad\quad $\blacksquare$\\
\end{dem}
%
%
%
%
%
%\addcontentsline{toc}{subsection}{1.1.3 The energy of a natural
%family of maps.}
%
\subsection{The energy of a natural family of maps.}
Let $a =(\theta, \cdots \!, 0)$  be a point of $\mathbf{B}^n$ with
$0<\nolinebreak[4]\theta<\nolinebreak[4]1$ and consider the map,
\[
u_a(x)= a+ \lambda_a(x)(x-a),
\]
where $\lambda_a(x)>0$ is chosen so that $u_a(x) \in
\mathbb{S}^{n-1}$ for any  $x \in \mathbf{B}^n \setminus  \{0\}$,
\[
\lambda_a(x)= \frac{\sqrt{\Delta_a(x)}-(a|x-a)}{\|x-a\|^2}
\]
and
\[
\Delta_a(x)=(1-\|a\|^2)\|x-a\|^2+ (a|x-a)^2.
\]
Notice that $u_a(x)=x$ as soon as $x$ is on the sphere. If we
denote by $\{e_i\}_{i \in \{1,\cdots,n\}}$ the standard basis of
$\mathbb{R}^n$, then,  $\forall i \leq n$, we have,
\begin{eqnarray*}
\|du_a(x).e_i\|^2 &=&
\left(\frac{\sqrt{\Delta_a}-(a|x-a)}{\|x-a\|^2}\right)^2\\
&+&
\Bigg[-2
\frac{(x-a|e_i)}{\|x-a\|^4}\Big(\sqrt{\Delta_a}-(a|x-a)\Big)\\
&&
+ \quad
\frac{(1-\|a\|^2)(x-a|e_i)+(x-a|a)(a|e_i)}{\sqrt{\Delta_a}\|x-a\|^2}\\
& &
-\frac{(a|e_i)}{\|x-a\|^2} \Bigg]^2\|x-a\|^2 \\
&+&
2 \Bigg(\frac{\sqrt{\Delta_a}-(a|x-a)}{\|x-a\|^2}\Bigg) \Bigg(-2
\frac{(x-a|e_i)}{\|x-a\|^4}\Big(\sqrt{\Delta_a}-(a|x-a)\Big)\\
&&
+
\quad\frac{(1-\|a\|^2)(x-a|e_i)+(x-a|a)(a|e_i)}{\sqrt{\Delta_a}\|x-a\|^2}\\
&&
-\frac{(a|e_i)}{\|x-a\|^2}\Bigg) (x-a|e_i).
\end{eqnarray*}
Let us prove that, for each $\alpha \in (-n,0)$,
$\int_{\mathbf{B}^n}r^\alpha \|\nabla u_a\|dx$ is finite. Consider
the map:
\begin{eqnarray*}
F: &\mathbb{R}^+ \times \mathbb{S}^{n-1}& \longrightarrow
\mathbb{R}^n\\ &(r,s)& \longmapsto a+rs=x.
\end{eqnarray*}
Then, we have,
\[
F^* \big(\|\bigtriangledown u_a\|^2 dx\big)=
\frac{1}{r^2}\dmath\sum_{i=1}^n H_{i,a}(s)\, r^{n-1}dr \wedge ds,
\]
where $H_{i,a}(s)$ is given on the sphere by,
\begin{eqnarray*}
H_{i,a}(s) &=& \left((1-\|a\|^2+ (a|s)^2) ^{1/2} - (a|s)
\right)^2\\ &+&  \bigg[ -2(s|e_i)\left((1-\|a\|^2+ (a|s)^2)
^{1/2}-(s|a)\right)\\ & &  +\quad \frac{(1-\|a\|^2)(s|e_i)+
(a|e_i)(s|a)}{(1-\|a\|^2+ (a|s)^2) ^{1/2}} - (a|e_i) \bigg]^2 \\
&+& 2 \bigg((1-\|a\|^2+ (a|s)^2) ^{1/2} - (a|s) \bigg)\\
 & & \bigg(-2(s|e_i)\left((1-\|a\|^2+ (a|s)^2) ^{1/2}-(s|a)\right)\\
 & &+\quad \frac{(1-\|a\|^2)(s|e_i)
  +
(a|e_i)(s|a)}{(1-\|a\|^2+ (a|s)^2) ^{1/2}} - (a|e_i)\bigg)\,
(s|e_i).
\end{eqnarray*}
It is clear that  $ H_{i,a}(s)$ is continuous on
$\mathbb{S}^{n-1}$. Therefore, near the point $a$, as $n \geq 3$,
the map $\|x\|^\alpha \|\nabla u_a\|$ is integrable. Furthermore,
near the point 0, as $\!\alpha>\!-n$, this map is also integrable.
In conclusion, for any $\alpha \in (-n,0)$, the energy
$E_{r^{\alpha},2}(u_a)$ is finite.\\
\newline
\textbf{Proof of Theorem 1.3(i)}. Since we have
\[
E_{2,r^{\alpha}}(u_0)=\int_{\mathbf{B}^n}\|x\|^\alpha \|\nabla
u_0\|^2 dx= \frac{|\mathbb{S}^{n-1}|(n-1)}{n+\alpha-2},
\]
the energy $E_{2,r^{\alpha}}(u_0)$ goes to infinity as $\alpha
\rightarrow -n+2$. On the other hand, as the energy
$E_{2,r^{\alpha}}(u_a)$ is continuous in $\alpha$, there exists a
real number $\alpha_0 \in (-n+2,0)$ such that, $\forall \alpha, \,
2\!-\!n\!<\!\alpha \leq\alpha_0$,
\[
\int_{\mathbf{B}^n}\|x\|^\alpha \|\nabla u_0\|^2 dx >
\int_{\mathbf{B}^n}\|x\|^\alpha \|\nabla u_a\|^2 dx.
\]
\textbf{Proof of Theorem 1.3(ii)}. Since $a=(
\theta,0,\cdots,\,0)$, we will study the function,
\[
G(\theta)=E_{2,r^{\alpha}}(u_a)=\int_{\mathbf{B}^n}r^{\alpha}
\|\nabla u_a\|^2 dx.
\]
Precisely, we will show that for any $\alpha \in (5-n,4-n)$, $G$
is two times differentiable at $\theta=0$ with
$\frac{dG}{d\theta}(0)=0$ and, when $\alpha$ is sufficiently close
to $4-n$, $\frac{d^2G}{d\theta^2}(0)<0$. Assertion (ii) of Theorem
1.3 then follows immediately. We have,
\begin{eqnarray*}
H_{i,a}(s)=H_{i,\theta}(s) &=& \bigg( \sqrt{
1\!-\!\theta^2\!+\!\theta^2s_1^2}-\theta s_1 \bigg)^2\\ &+&
\bigg(\!-\!2s_i\big(\sqrt{
1\!-\!\theta^2\!+\!\theta^2s_1^2}-\theta s_1\big) +
\frac{(1-\theta^2)s_i+\delta_{i1}\theta^2s_1}{\sqrt{
1\!-\!\theta^2\!+\!\theta^2s_1^2}}-\delta_{i1}\theta \bigg)^2\\
&+& 2\bigg(\sqrt{ 1\!-\!\theta^2\!+\!\theta^2s_1^2}-\theta s_1
\bigg)\bigg(\!-\!2s_i\big(\sqrt{
1\!-\!\theta^2\!+\!\theta^2s_1^2}-\theta s_1 \big)\\
&&
\hspace{3cm}
+\,\frac{(1-\theta^2)s_i+\delta_{i1}\theta^2s_1}{\sqrt{
1\!-\!\theta^2\!+\!\theta^2s_1^2}}-\delta_{i1}\theta \bigg)s_i,
\end{eqnarray*}
where $\delta_{ij}= 0$ if $i \neq j$ and $0$ else.\\
We notice that $H_{i,\theta}(s)$ is bounded on $[0,1] \times
\mathbb{S}^{n-1}$. Indeed, for all $x,y,z \in [0,1]$, excepting
$(x,y)=(0,1)$, we have,
\[
\left|\frac{x}{\sqrt{1-y^2+y^2x^2)}}\right|\leq 1 \quad \text{and}
\quad \left|\frac{(1-y^2)z}{\sqrt{1-y^2+y^2x^2)}}\right|\leq 1 .
\]
Then, for almost all $(s,\theta) \in \mathbb{S}^{n-1} \times
[0,1]$, we have,
\[
\left|\frac{(1-\theta^2)s_i+\delta_{i1}\theta^2s_1}{\sqrt{
1\!-\!\theta^2\!+\!\theta^2s_1^2}}\right|\leq 1,
\]
and the others terms are continuous in $ [0,1] \times
\mathbb{S}^{n-1}$.\\
We have,
\begin{eqnarray*}
E_{2,r^\alpha}(u_a)&\!=\!&\!\int_{\mathbf{B}^n}\!\|x\|^\alpha
\|\nabla u_a\|^2 dx\!
=
\!\!\int_{\mathbf{B}^n}\!\!\|a+rs\|^{\alpha} r^{n-3}
H(\theta,s)drds\\
&=&
\int_{\mathbf{S}^{n-1}}\!H(\theta,s)\!\left(
\int_0^{\gamma_{\theta}(s)}\!\!\big((r+\theta
s_1)^2+\theta^2(1-s_1^2)\big)^{\alpha /2}r^{n-3}dr\!\! \right)\!
ds ,
\end{eqnarray*}
where $\gamma_{\theta}(s)=\sqrt{
1\!-\!\theta^2\!+\!\theta^2s_1^2}-\theta s_1$ and
$H(\theta,s)=\dmath\sum_{i=1}^n H_{i,\theta}(s)$.
We notice that $H(\theta,s)$ is indefinitely differentiable in
$(-1/2,1/2) \times \mathbb{S}^{n-1}$. Let $C_n$ be a positive real
number so that, $\forall (\theta,s) \in (-1/2,1/2) \times
\mathbb{S}^{n-1}$
\[
\left|H (\theta,s)\right| \leq C_n,\,\left|\frac{\partial
H(\theta,s)}{\partial \theta}\right| \leq C_n, \,
\left|\frac{\partial^2 H(\theta,s)}{\partial \theta^2}\right| \leq
C_n.
\]
Furthermore, we have,
\[
H(\theta,s)=(n\!-\!1) -2(n\!-\!1)s_1\theta + ((2n\!-\!3)s_1^2
\!-n+2)\theta^2 +o(\theta^2).\quad (A)
\]
Let us set $\rho=r+\theta s_1$, $\beta(\theta,s)=\sqrt{
1\!-\!\theta^2\!+\!\theta^2s_1^2}$ and
\[
F(\theta,s)= \int_{\theta s_1}^{\beta(\theta,s)} (\rho-\theta
s_1)^{n-3} (\rho^2+\theta^2(1-s_1^2))^{\alpha/2}d\rho.
\]
Notice that $\rho \in [-1,3]$. Then, $G(\theta)= \int_{
\mathbb{S}^{n-1}}H(\theta,s)F(\theta,s) ds$.
Let us set $g(\rho,\theta,s)=(\rho-\theta s_1)^{n-3}
(\rho^2+\theta^2(1-s_1^2))^{\alpha/2}$.
\begin{lem}
The map $\theta \mapsto G(\theta)$ is continuous on $(-1/2,1/2)$
and continuously differentiable on $(-1/2,1/2)\!\setminus \!\{0\}$
for any $\alpha > 3\!-\!n$ .
\end{lem}
\begin{dem}
We have, $\forall s \in \mathbb{S}^{n-1}\!\setminus\!\{(\pm
1,0,\cdots,0)\}$,
\[
\frac{(\rho-\theta s_1)^{2}}{(\rho^2+\theta^2(1-s_1^2))}\leq
\frac{2}{1-s_1^2} \quad (1.1)
\]
Indeed, $(1-s_1^2)(\rho-\theta s_1)^2\leq
2(1-s_1^2)(\rho^2+\theta^2) \leq 2(\rho^2+\theta^2(1-s_1^2))$. And
then,
\[
g(\rho,\theta,s)\leq
\frac{2^{\frac{n-3}{2}}}{(1-s_1^2)^{\frac{n-3}{2}}}(\rho^2+\theta^2(1-s_1^2))^{\frac
{\alpha+n-3}{2}}. \quad (1.2)
\]
Since $\alpha >3-n$ we deduce that the map
$(\rho,\theta)\rightarrow g(\rho,\theta,s)$ is continuous on
$(-1/2,1/2) \times [-1,3]$. Hence, the map $z \mapsto \int_0^z g(
\rho,\theta,s)d\rho$ is differentiable on $[-1,3]$ and,
\[
\frac{\partial}{\partial z}\int_0^z g( \rho,\theta,s)d\rho=
g(z,\theta,s).
\]
Furthermore, for any $\rho \in [-1,3]$, the map $\theta \mapsto
g(\rho,\theta,s)$ is differentiable and
\begin{eqnarray*}
\frac{\partial g}{\partial \theta }(\rho,\theta,s)&=&
-(n-3)s_1(\rho-\theta
s_1)^{n-4}(\rho^2+\theta^2(1-s_1^2))^{\frac{\alpha}{2}}
\\&+&
\frac{\alpha}{2}(\rho-\theta
s_1)^{n-3}2\theta(1-s_1^2)(\rho^2+\theta^2(1-s_1^2))^{\frac{\alpha}{2}-1}.
\end{eqnarray*}
Let $a,b$ be two real in $(0,1/2)$ with $a<b$. We have for any
$|\theta| \in (a,b)$, for any $ s \in
\mathbb{S}^{n-1}\!\setminus\!\{(\pm 1,0,\cdots,0)\}$,
\begin{eqnarray*}
\left|\frac{\partial g}{\partial \theta
}(\rho,\theta,s)\right|&\leq&
(n-3)4^{n-4}(a^2(1-s_1^2))^{\frac{\alpha}{2}}
\\&+&
|\alpha|4^{n-3}(1-s_1^2)(a^2(1-s_1^2))^{\frac{\alpha}{2}-1}.\quad
(1.3)
\end{eqnarray*}
This shows that $\theta \mapsto \int_0^z g(\rho,\theta,s)d\rho$ is
differentiable on $(-1/2,1/2)\setminus \{0\}$ and
\[
\frac{\partial}{\partial \theta}\int_0^z g(\rho,\theta,s)d\rho=
\int_0^z \frac{\partial g}{\partial \theta }(\rho,\theta,s)d\rho.
\]
Moreover the map $(z,\theta)\mapsto  \int_0^z \frac{\partial g
}{\partial \theta}(\rho,\theta,s)d\rho$ is continuous in
$[-1,3]\times (-1/2,1/2)\!\setminus \!\{0\}$. Indeed, $\theta
\mapsto \frac{\partial g}{\partial \theta } (\rho,\theta,s)$ is
clearly continuous on $(-1/2,1/2)\!\setminus \!\{0\}$ and from
(1.3) and by Lebesgue Theorem,  $\theta \mapsto \int_0^z
\frac{\partial g}{\partial \theta }(\rho,\theta,s)d\rho$ is
continuous on $(-1/2,1/2)\setminus \{0\}$. Then, for any $\epsilon
>0$, we will have for any
sufficiently small $h,k$,
\begin{eqnarray*}
\left| \int_0^{z+h} \frac{\partial g}{\partial \theta
}(\rho,\theta+k,s)d\rho -\int_0^z \frac{\partial g}{\partial
\theta }(\rho,\theta,s)d\rho \right|
&\leq&
\Big| \int_0^{z} \frac{\partial g}{\partial \theta}
(\rho,\theta\!+\!k,s)d\rho
\\
&&
- \int_0^z \frac{\partial g}{\partial \theta }(\rho,\theta,s)d\rho
\Big|\\
&+&
\left|\int_z^{z+h}\!\! \frac{\partial g}{\partial \theta
}(\rho,\theta\!+\!k,s)d\rho \right|\\
&\leq&
\epsilon.
\end{eqnarray*}
The map $(z,\theta) \mapsto \int_0^z g(\rho,\theta,s)d\rho$ is
differentiable on $[-1,3]\times (-1/2,1/2)\!\setminus \!\{0\}$ and
the map $\theta \mapsto F(\theta,s)$ is differentiable in
$(-1/2,1/2)\!\setminus\! \{0\}$ and for any $\theta \in
(-1/2,1/2)\!\setminus \!\{0\}$,
\begin{eqnarray*}
\frac{\partial F}{\partial \theta }(\theta,s)
&=&
\frac{\partial \beta}{\partial \theta
}(\theta,s)g(\beta(\theta,s),\theta,s)-s_1g(\theta
s_1,\theta,s)+\int_{\theta s_1}^{\beta(\theta,s)}\frac{\partial g
}{\partial \theta }(\rho,\theta,s)d\rho \\
&=&
\frac{\theta(s_1^2-1)}{(1-\theta^2+\theta^2s_1^2)^{1/2}}((1-\theta^2+\theta^2s_1^2)^{1/2}-\theta
s_1)^{n-3}\\
&+&
\int_{\theta s_1}^{\beta(\theta,s)}g_1(\rho,\theta,s)d\rho
+
\int_{\theta s_1}^{\beta(\theta,s)}g_2(\rho,\theta,s)d\rho,
\end{eqnarray*}
where,
\[
g_1(\rho,\theta,s)=-(n-3)s_1(\rho-\theta
s_1)^{n-4}(\rho^2+\theta^2(1-s_1^2))^{\frac{\alpha}{2}}
\]
and
\[
g_2(\rho,\theta,s)=\frac{\alpha}{2}(\rho-\theta
s_1)^{n-3}2\theta(1-s_1^2)(\rho^2+\theta^2(1-s_1^2))^{\frac{\alpha}{2}-1}.
\]
Now, the map $\theta \mapsto F(\theta,s)$ is continuous on
$(-1/2,1/2)$. Indeed, since the map  $ \theta \mapsto
g(\rho,\theta,s)d\rho$ is continuous on $(-1/2,1/2)$ and from
(1.2) $ \theta \mapsto \int_0^zg(\rho,\theta,s)d\rho$ is
continuous on $(-1/2,1/2)$. Then, for any $\epsilon > 0$, we have
$ \forall h,k$ sufficiently small,
\begin{eqnarray*}
\left| \int_0^{z+h} g(\rho,\theta+k,s)d\rho- \int_0^z
g(\rho,\theta,s)d\rho \right|
&\leq&
\Big| \int_0^{z} g(\rho,\theta\!+\!k,s)d\rho\\
&&
-\int_0^z g(\rho,\theta,s)d\rho \Big|\\
&+&
\left|\int_z^{z+h} g(\rho,\theta\!+\!k,s)d\rho \right|\\
&\leq&
\epsilon.
\end{eqnarray*}
Then, the map $(z,\theta)\mapsto \int_0^z g(\rho,\theta,s)d\rho$
is continuous on $[-1,3] \times (-1/2,1/2)$ and consequently
$\theta \mapsto F(\theta,s)$ is continuous on $(-1/2,1/2)$.

Now, we know that $\theta \mapsto H(\theta,s)F(\theta,s)$ is
continuous on $(-1/2,1/2)$ and differentiable on
$(-1/2,1/2)\setminus \{0\}$. Furthermore from (1.2), we have, for
any $s \in \mathbb{S}^{n-1}\!\setminus\!\{( \pm 1,0,\cdots,0)\}$,
\[
\left|H(\theta,s)F(\theta,s) \right| \leq
3.2^{\frac{n-3}{2}}10^{\frac{\alpha+n-3}{2}}C_n.\frac{1}{(1-s_1^2)^{\frac{n-3}{2}}}.\quad
(1.4)
\]
\[
\left|\frac{\partial H}{\partial \theta }(\theta,s)F(\theta,s)
\right| \leq
3.2^{\frac{n-3}{2}}10^{\frac{\alpha+n-3}{2}}C_n.\frac{1}{(1-s_1^2)^{\frac{n-3}{2}}}.
\quad(1.5)
\]
Consider the map $\eta:(\theta,s)\mapsto \eta(\theta,s)=
\frac{\theta(s_1^2-1)}{\sqrt{1-\theta^2+\theta^2s_1^2}}((1-\theta^2+\theta^2s_1^2)^{1/2}-\theta
s_1)^{n-3}$. This map is indefinitely differentiable on
$(-1/2,1/2) \times \mathbb{S}^{n-1}$. Let $B_{n}$ be a positive
real number so that, $\forall (\theta,s) \in (-1/2,1/2) \times
\mathbb{S}^{n-1}$,
\[
\left| \eta(\theta,s) \right| \leq B_n \quad \left| \frac{\partial
\eta}{\partial \theta }(\theta,s) \right| \leq B_n.
\]
Considering $a,b \in (0,1/2)$ with $a<b$ we have, for any $\theta
\in (a,b)$, for any $s \in \mathbb{S}^{n-1}\!\setminus\!\{(\pm
1,0,\cdots,0)\}$,
\begin{eqnarray*}
\left|H(\theta,s)\frac{\partial F}{\partial \theta }(\theta,s)
\right|
& \leq&
\Big(B_n +3(n-3).4^{n-4}.a^{\alpha}(1-s_1^2)^{\frac{\alpha}{2}}\\
&+&
\!|3\alpha|.4^{n-3}a^{\alpha-1}(1-s_1^2)^{\frac{\alpha}{2}}\Big)
C_n.\quad (1.6)
\end{eqnarray*}
Since the maps $s \mapsto \frac{1}{(1-s_1^2)^{\frac{n-3}{2}}}$ and
$s \mapsto (1-s_1^2)^{\frac{\alpha}{2}}$ are integrable on
$\mathbb{S}^{n-1}$, we deduce that $\theta \mapsto G(\theta)$ is
continuous on $(-1/2,1/2)$ and continuously differentiable  on
$(-1/2,1/2)\setminus \{0\}$.
\end{dem}
\begin{lem}
The map $\theta \mapsto G(\theta)$ is differentiable  at $0$ and
$\frac{d G}{d \theta } (0)= 0$.
\end{lem}
\begin{dem}
Since for any $s \in \mathbb{S}^{n-1}\!\setminus\! \{(\pm
1,0,\cdots,0)\}$, $\theta \mapsto F(\theta,s)$ is continuous on
$(-1/2,1/2)$ from $(A)$ we have,
\[
\frac{\partial H}{\partial \theta }(\theta,s)F(\theta,s)
\mathop{\longrightarrow}_{\theta \rightarrow 0} \frac{\partial H
}{\partial \theta
}(0,s)F(0,s)=-2(n-1)s_1\int_0^1\rho^{n-3+\alpha}d\rho=\frac{-2(n-1)s_1}{n-2+\alpha}.
\]
From $(1.5)$ and Lebesgue Theorem we have,
\[
\int_{\mathbb{S}^{n-1}}\frac{\partial H}{\partial \theta
}(\theta,s)F(\theta,s)ds\mathop{\longrightarrow}_{\theta
\rightarrow 0}
\int_{\mathbb{S}^{n-1}}\frac{-2(n-1)s_1}{n-2+\alpha}ds=0.
\]
Moreover, it is clear that,
\[
\int_{\mathbb{S}^{n-1}}H(\theta,s)\eta(\theta,s)ds\mathop{\longrightarrow}_{\theta
\rightarrow 0} 0.
\]
Let $J(m,n)$ be the integral,
\[
J(m,n)=\int_{\frac{s_1}{\sqrt{1-s_1^2}}}^{\sqrt{\frac{1}{\theta^2(1-s_1^2)}-1}}
(\sqrt{1-s_1^2}\,t-s_1)^{m}(t^2+1)^{n}dt.
\]
Notice that $J(m,n)$ converges as $\theta$ goes to $0$ if and only
if $m+2n< -1$.
Consider the change of variables $\rho=t\theta \sqrt{1-s_1^2}$ if
$\theta>0$. If $\theta <0$, then we set $\rho=-t\theta
\sqrt{1-s_1^2}$ and conclusion will be the same. Hence, we assume
that $\theta >0$. Then,
\begin{eqnarray*}
\int_{\theta s_1}^{\beta(\theta,s)}g_1(\rho,\theta,s) d\rho
&=&
-(n-3)s_1(1-s_1^2)^{\frac{1+\alpha}{2}}
\theta^{n-3+\alpha}J(n-4,\frac{\alpha}{2}).
\end{eqnarray*}
\begin{eqnarray*}
\int_{\theta s_1}^{\beta(\theta,s)}g_2(\rho,\theta,s) d\rho
&=&
\alpha\theta^{n-3+\alpha}(1-s_1^2)^{\frac{1+\alpha}{2}}J(n-3,\frac{\alpha}{2}-1).
\end{eqnarray*}
\textbf{First case:\,$\alpha \geq 4-n$}.

$J(n-4,\frac{\alpha}{2})$ and $J(n-3,\frac{\alpha}{2}-1)$ go to $+
\infty$ as $\theta \rightarrow 0$. Furthermore, we have,
\[
J(n-4,\frac{\alpha}{2})
\mathop{\sim}_{0}
(1-s_1^2)^{\frac{n-4}{2}}\int_{\frac{s_1}{\sqrt{1-s_1^2}}}^{\sqrt{\frac{1}{\theta^2(1-s_1^2)}-1}}t^{n-4+\alpha}dt
\]
\[
J(n-4,\frac{\alpha}{2})
\mathop{\sim}_{0}
\frac{1}{n-3+\alpha}\frac{1}{\theta^{n-3+\alpha}}(1-s_1^2)^{\frac{-1-\alpha}{2}}.
\]
Since $t^{n+\alpha-5}$ may be equal to zero at zero, we write,
\begin{eqnarray*}
J(n-3,\frac{\alpha}{2}-1)
&=&
\int_{\frac{s_1}{\sqrt{1-s_1^2}}}^1(\sqrt{1-s_1^2}\,t-s_1)^{n-3}(t^2+1)^{\frac{\alpha}{2}-1}dt\\
&&
+\int_1^{\sqrt{\frac{1}{\theta^2(1-s_1^2)}-1}}
(\sqrt{1-s_1^2}\,t-s_1)^{n-3}(t^2+1)^{\frac{\alpha}{2}-1}dt.
\end{eqnarray*}
We have,
\[
J(n-3,\frac{\alpha}{2}-1)
\mathop{\sim}_{0}
(1-s_1^2)^{\frac{n-3}{2}}\int_1^{\sqrt{\frac{1}{\theta^2(1-s_1^2)}-1}}t^{n-5+\alpha}dt.
\]
Then, if $\alpha \neq 4-n$,
\[
J(n-3,\frac{\alpha}{2}-1)
\mathop{\sim}_{0}
\frac{1}{n-4+\alpha}\frac{1}{\theta^{n-4+\alpha}}(1-s_1^2)^{\frac{1-\alpha}{2}},
\]
and note that if $\alpha=4-n$,
$J(n-3,\frac{\alpha}{2}-1)
\mathop{\sim}_{0}
-(1-s_1^2)^{\frac{n-3}{2}}\ln(\theta^2(1-s_1^2)).$
Hence, by $(A)$ we have,
\begin{eqnarray*}
H(\theta,s)\int_{\theta
s_1}^{\beta(\theta,s)}g_1(\rho,\theta,s)d\rho
&=&
-H(\theta,s)(n-3)s_1(1-s_1^2)^{\frac{\alpha+1}{2}}
\theta^{n-3+\alpha}I_1\\
&&
\mathop{\longrightarrow}_{\theta \rightarrow 0}
-\frac{(n-3)(n-1)}{n-3+\alpha}s_1,
\end{eqnarray*}

and
\[
H(\theta,s)\int_{\theta
s_1}^{\beta(\theta,s)}g_2(\rho,\theta,s)d\rho=H(\theta,s)\alpha(1-s_1^2)^{\frac{\alpha+1}{2}}
\theta^{n-3+\alpha}I_2 \mathop{\longrightarrow}_{\theta
\rightarrow 0} 0.
\]
Observe that $ \frac{|s_1|}{\sqrt{1-s_1^2}}\leq
\sqrt{\frac{1}{\theta^2(1-s_1^2)}-1}$. Indeed, $s_1^2\theta^2 \leq
1-\theta^2+\theta^2s_1^2$. It follows from $(1.1)$ that
\[
(\rho- \theta s_1)^{n-4}\leq
\frac{2^{\frac{n-4}{2}}}{(1-s_1^2)^{\frac{n-4}{2}}}(\rho^2+
\theta^2(1-s_1^2))^{\frac{n-4}{2}}.
\]
Recall that $\rho=t\theta \sqrt{1-s_1^2}$. Since $\alpha \geq 4-n$
, we have, for any $s \in \mathbb{S}^{n-1}\!\setminus\!\{(\pm
1,0,\cdots,0)\}$,
\begin{eqnarray*}
\left|H(\theta,s)\int_{\theta
s_1}^{\beta(\theta,s)}g_1(\rho,\theta,s)d \rho \right|
&\leq&
2C_n(n\!-\!3)\!2^{\frac{n-4}{2}}(1-s_1^2)^{\frac{\alpha+1}{2}}
\theta^{n-3+\alpha}\\
&&
\times\int_0^{\sqrt{\frac{1}{\theta^2(1-s_1^2)}-1}}(t^2+1)^{\frac{n-4+\alpha}{2}}dt\\
&\leq&
C_n(n\!-\!3)2^{\frac{n-2}{2}}(1\!-\!s_1^2)^{\frac{\alpha+1}{2}}\theta^{n-3+\alpha}\\
&&
\times\sqrt{\frac{1}{\theta^2(1\!-\!s_1^2)}\!\!-\!\!1}
\!\!\!\left(\!\frac{1}{\theta^2(1-s_1^2)}\!\right)^{\frac{n-4+\alpha}{2}}\\
&\leq&
C_n(n-3)2^{\frac{n-2}{2}}(1-s_1^2)^{\frac{-n+4}{2}}\sqrt{1-\theta^2(1-s_1^2)}\\
&\leq&
C_n(n-3)2^{\frac{n-1}{2}}(1-s_1^2)^{\frac{-n+4}{2}}.
\end{eqnarray*}
Since $s\mapsto (1-s_1^2)^{\frac{-n+4}{2}}$ is integrable on
$\mathbb{S}^{n-1}$, by Lebesgue Theorem we have,
\[
\int_{\mathbb{S}^{n-1}}H(\theta,s)\int_{\theta
s_1}^{\beta(\theta,s)}g_1(\rho,\theta,s)d\rho ds
\mathop{\longrightarrow}_{\theta \rightarrow
0}-\int_{\mathbb{S}^{n-1}}\frac{(n-3)(n-1)}{n-3+\alpha}s_1 ds = 0.
\]
Moreover, we have, for any $s \in
\mathbb{S}^{n-1}\!\setminus\!\{(\pm 1,0,\cdots,0)\}$, since
$\alpha +n-5 \geq 0$,
\begin{eqnarray*}
\left|H(\theta,s)\int_{\theta
s_1}^{\beta(\theta,s)}g_2(\rho,\theta,s)d \rho \right|
&\leq&
2C_n|\alpha|2^{\frac{n-3}{2}}(1-s_1^2)^{\frac{\alpha+1}{2}}
\theta^{n-3+\alpha}\\
&&
\times\int_0^{\sqrt{\frac{1}{\theta^2(1-s_1^2)}-1}}(t^2+1)^{\frac{n-5+\alpha}{2}}dt\\
&\leq&
C_n|\alpha|2^{\frac{n-2}{2}}(1\!-\!s_1^2)^{\frac{\alpha+1}{2}}\theta^{n-3+\alpha}
\!\!\!\!\\
&&
\times\int_0^{\sqrt{\frac{1}{\theta^2(1-s_1^2)}-1}}\!\!\!\frac{1}{(t^2+1)}dt
\left(\!\frac{1}{\theta^2(1-s_1^2)}\!\right)^{\frac{n-3+\alpha}{2}}\\
&\leq&
C_n|\alpha|2^{\frac{n-2}{2}}\frac{\pi}{2}(1-s_1^2)^{\frac{-n+4}{2}}.
\end{eqnarray*}
Then, by Lebesgue Theorem,
\[
\int_{\mathbb{S}^{n-1}}H(\theta,s)\int_{\theta
s_1}^{\beta(\theta,s)}g_2(\rho,\theta,s)d\rho ds
\mathop{\longrightarrow}_{\theta \rightarrow 0} 0.
\]

\textbf{Second case:\,$3-n<\alpha < 4-n$}.

For the same reasons that when $\alpha \geq 4-n$, we have,
\[
H(\theta,s)\int_{\theta
s_1}^{\beta(\theta,s)}g_1(\rho,\theta,s)d\rho
\mathop{\longrightarrow}_{\theta \rightarrow 0}
-\frac{(n-3)(n-1)}{n-3+\alpha}s_1.
\]
Furthermore, as $4-n> \alpha > 3-n$, $\forall s \in
\mathbb{S}^{n-1}\!\setminus\!\{(-1,0,\cdots,0),(1,0,\cdots,0)\}$,
\begin{eqnarray*}
\left|H(\theta,s)\int_{\theta
s_1}^{\beta(\theta,s)}g_1(\rho,\theta,s)d \rho \right|
&\leq&
2C_n(n-3)2^{\frac{n-4}{2}}(1-s_1^2)^{\frac{\alpha+1}{2}}
\theta^{n-3+\alpha}\\
&&
\times\int_0^{\sqrt{\frac{1}{\theta^2(1-s_1^2)}-1}}
(t^2+1)^{\frac{n-4+\alpha}{2}}dt\\
&\leq&
C_n(n-3)2^{\frac{n-2}{2}}(1\!-\!s_1^2)^{\frac{\alpha+1}{2}}\theta^{n-3+\alpha}\\
&&
\times\int_0^{\sqrt{\frac{1}{\theta^2(1-s_1^2)}-1}}(t^2)^{\frac{n-4+\alpha}{2}}dt
\\
&\leq&
\frac{C_n(n-3)2^{\frac{n-2}{2}}(1-s_1^2)^{\frac{\alpha+1}{2}}\theta^{n-3+\alpha}}{n-3+
\alpha} \big(
\frac{1}{\theta^2(1-s_1^2)}-1\big)^{\frac{n-3+\alpha}{2}}.\\
&\leq&
\frac{C_n(n-3)2^{\frac{2n-7+\alpha}{2}}(1-s_1^2)^{\frac{4-n}{2}}}{n-3+
\alpha}.
\end{eqnarray*}
Then, by Lebesgue Theorem,
\[
\int_{\mathbb{S}^{n-1}}H(\theta,s)\int_{\theta
s_1}^{\beta(\theta,s)}g_1(\rho,\theta,s)d\rho ds
\mathop{\longrightarrow}_{\theta \rightarrow
0}-\int_{\mathbb{S}^{n-1}}\frac{(n-3)(n-1)}{n-3+\alpha}s_1 ds=  0.
\]
Moreover, $J(n-3,\frac{\alpha}{2}-1)$ is finite when $\theta
\rightarrow 0$ then, as $\alpha > 3-n$,
Furthermore,
\[
H(\theta,s)\int_{\theta
s_1}^{\beta(\theta,s)}g_2(\rho,\theta,s)d\rho ds
\mathop{\longrightarrow}_{\theta \rightarrow 0}  0.
\]
\begin{eqnarray*}
\left|H(\theta,s)\int_{\theta
s_1}^{\beta(\theta,s)}g_2(\rho,\theta,s)d \rho \right|
&\leq&
2C_n|\alpha|2^{\frac{n-3}{2}}(1-s_1^2)^{\frac{\alpha+1}{2}}
\theta^{n-3+\alpha}\\
&&
\times\int_0^{\sqrt{\frac{1}{\theta^2(1-s_1^2)}-1}}(t^2+1)^{\frac{n+\alpha-5}{2}}dt\\
&\leq&
C_n|\alpha|2^{\frac{n-1}{2}}(1\!-\!s_1^2)^{\frac{\alpha+1}{2}}\theta^{n-3+\alpha}\\
&&
\times\int_0^{\sqrt{\frac{1}{\theta^2(1-s_1^2)}-1}}\!\!\!\frac{1}{(t^2+1)}dt
\left(\!\frac{1}{\theta^2(1-s_1^2)}\!\right)^{\frac{n-3+\alpha}{2}}\\
&\leq&
C_n|\alpha|2^{\frac{n-1}{2}}(1-s_1^2)^{\frac{-n+4}{2}}\int_0^{+\infty}\frac{1}{(t^2+1)}dt.
\end{eqnarray*}
Then, by Lebesgue Theorem,
\[
\int_{\mathbb{S}^{n-1}}H(\theta,s)\int_{\theta
s_1}^{\beta(\theta,s)}g_2(\rho,\theta,s)d\rho ds
\mathop{\longrightarrow}_{\theta \rightarrow 0} 0.
\]
Finally, we have
\[
\frac{d G}{d \theta
 }(\theta)\mathop{\longrightarrow}_{\theta \rightarrow 0} 0. %
\]
By Lemma 1.4 we deduce that $G$ is differentiable at $0$ and
$\frac{d G}{d \theta }(0)=0$.
\end{dem}
\begin{lem}
The map $\theta \rightarrow G(\theta)$ is two times differentiable
on $(-1/2,1/2)\!\setminus\!\{0\}$.
\end{lem}
\begin{dem}
We know that the map $\theta \rightarrow \frac{\partial
H}{\partial \theta }(\theta,s)F(\theta,s)$ is differentiable on
$(-1/2,1/2)\!\setminus\!\{0\}$.
The maps $\theta \rightarrow \eta(\theta,s)$, $\theta \rightarrow
g_1(\rho,\theta,s)$, $\theta \rightarrow g_2(\rho,\theta,s)$ are
differentiable on $(-1/2,1/2)\!\setminus\!\{0\}$. We have,
\begin{eqnarray*}
\frac{\partial \eta}{\partial \theta }(\theta,s)
&=&
\frac{(s_1^2-1)\sqrt{1-\theta^2+\theta^2s_1^2}-\theta(s_1^2-1)
\frac{\theta(s_1^2-1)}{\sqrt{1-\theta^2+\theta^2s_1^2}}}
{1-\theta^2+\theta^2s_1^2}\\
&&
\hspace{3cm}\times(\sqrt{1-\theta^2+\theta^2s_1^2}-\theta
s_1)^{n-3}\\
&&
+\,\frac{(n-3)\theta(s_1^2-1)}{\sqrt{1-\theta^2+\theta^2s_1^2}}
\left(\frac{\theta(s_1^2-1)}{\sqrt{1-\theta^2+\theta^2s_1^2}}
-s_1\right)\\
&&
\hspace{3cm}\times(\sqrt{1-\theta^2+\theta^2s_1^2}-\theta
s_1)^{n-4}.
\end{eqnarray*}
\begin{eqnarray*}
\frac{\partial g_1}{\partial \theta} (\rho,\theta,s)
&=&
(n-3)(n-4)s_1^2(\rho-\theta
s_1)^{n-5}(\rho^2+\theta^2(1-s_1^2))^{\frac{\alpha}{2}}\\
&&
-\, \alpha(n-3)s_1(1-s_1^2)\theta(\rho-\theta
s_1)^{n-4}(\rho^2+\theta^2(1-s_1^2))^{\frac{\alpha}{2}-1}.
\end{eqnarray*}
\begin{eqnarray*}
\frac{\partial g_2}{\partial \theta} (\rho,\theta,s)
&=&
-\alpha (n-3) s_1(1-s_1^2)\theta (\rho-\theta
s_1)^{n-4}(\rho^2+\theta^2(1-s_1^2))^{\frac{\alpha}{2}-1}\\
&&
\quad  +\, \alpha(\alpha-2)(1-s_1^2)^2\theta^2(\rho-\theta
s_1)^{n-3}(\rho^2+\theta^2(1-s_1^2))^{\frac{\alpha}{2}-2}\\
&&
\quad  + \,\alpha(1-s_1^2)(\rho-\theta
s_1)^{n-3}(\rho^2+\theta^2(1-s_1^2))^{\frac{\alpha}{2}-1}.
\end{eqnarray*}
We set,
\[
g_{11}(\rho,\theta,s)= (n-3)(n-4)s_1^2(\rho-\theta
s_1)^{n-5}(\rho^2+\theta^2(1-s_1^2))^{\frac{\alpha}{2}},
\]
\[
g_{12}(\rho,\theta,s)=-2\alpha(n-3)s_1(1-s_1^2)\theta(\rho-\theta
s_1)^{n-4}(\rho^2+\theta^2(1-s_1^2))^{\frac{\alpha}{2}-1}.
\]
\[
g_{21}(\rho,\theta,s)=
\alpha(\alpha-2)(1-s_1^2)^2\theta^2(\rho-\theta
s_1)^{n-3}(\rho^2+\theta^2(1-s_1^2))^{\frac{\alpha}{2}-2},
\]
\[
g_{22}(\rho,\theta,s)=\alpha(1-s_1^2)(\rho-\theta
s_1)^{n-3}(\rho^2+\theta^2(1-s_1^2))^{\frac{\alpha}{2}-1}.
\]
Let $a,b \in (0,1/2)$ with $a<b$. We have, $\forall s \in
\mathbb{S}^{n-1}\setminus\{(\pm 1,0,\cdots,0)\}$,
\[
\left|\frac{\partial g_1}{\partial \theta} (\rho,\theta,s)\right|
\leq (n-3)(n-4)4^{n-5}a^{\alpha}(1-s_1^2)^{\frac{\alpha}{2}}+
|\alpha|(n-3)4^{n-4}a^{\alpha-1}(1-s_1^2)^{\frac{\alpha}{2}}.\quad
(1.7)
\]
\begin{eqnarray*}
\left|\frac{\partial g_2}{\partial \theta} (\rho,\theta,s)\right|
&\leq&
|\alpha(\alpha-2)|4^{n-3}a^{\alpha-2}(1-s_1^2)^{\frac{\alpha}{2}}\\
&&
\quad  +\,
|\alpha|4^{n-3}a^{\alpha-1}(1-s_1^2)^{\frac{\alpha}{2}}\\
&&
\quad
+\,|\alpha|(n-3)4^{n-4}a^{\alpha-1}(1-s_1^2)^{\frac{\alpha}{2}}.\quad
(1.8)
\end{eqnarray*}
Then, for any $i \in \{1,2\}$, the maps $\theta \mapsto \int_0^z
g_i(\rho,\theta,s)d\rho$ is
 differentiable on $(0,1/2)$, and
\[
\frac{\partial }{\partial \theta}\int_0^z
g_i(\rho,\theta,s)d\rho=\int_0^z \frac{\partial g_i}{\partial
\theta}(\rho,\theta,s)d\rho.
\]
Furthermore, for any $i \in \{1,2\}$, $\theta \mapsto
\frac{\partial g_i}{\partial \theta}(\rho,\theta,s)$ is continuous
on $(-1/2,1/2)\!\setminus\!\{0\}$, then, $\theta \mapsto
\int_0^z\frac{\partial g_i}{\partial \theta}(\rho,\theta,s)d\rho$,
is continuous on $(-1/2,1/2)\!\setminus\!\{0\}$.
Hence, for any $i \in \{1,2\}$ and for any $\epsilon >0$, we have
$\forall h,k$ two sufficiently small,
\begin{eqnarray*}
\left| \int_0^{z+h} \frac{\partial g_i}{\partial \theta
}(\rho,\theta+k,s)d\rho- \int_0^z \frac{\partial g_i}{\partial
\theta }(\rho,\theta,s)d\rho \right|
&\leq&
\Big| \int_0^{z} \frac{\partial g_i}{\partial \theta }
(\rho,\theta\!+\!k,s)d\rho\\
&&
- \int_0^z \frac{\partial g_i}{\partial \theta
}(\rho,\theta,s)d\rho \Big|\\
&+&
\left|\int_z^{z+h}\!\! \frac{\partial g_i}{\partial \theta
}(\rho,\theta\!+\!k,s)d\rho \right|\\ &\leq& \epsilon.
\end{eqnarray*}
This proves that for any $i \in \{1,2\}$, $(z,\theta) \mapsto
\int_0^z\frac{\partial g_i}{\partial
\theta}(\rho,\theta,s)\!d\rho$ is continuous on $[-1,3]\times
(-1/2,1/2)\!\setminus\!\{0\}$.
Moreover,  for any $i \in \{1,2\}$ the map $\rho \mapsto
g_i(\rho,\theta,s)$ is continuous on $[-1,3]$ for any $\theta \in
(-1/2,1/2)\!\setminus\!\{0\}$. Then, $z \mapsto \int_0^z
g_i(\rho,\theta,s) d\rho$ is differentiable on $[-1,3]$ for any
$\theta  \in (-1/2,1/2)\!\setminus\!\{0\}$ and  $\frac{\partial
}{\partial z}\int_0^z g_i(\rho,\theta,s) d\rho=
g_i(z,\theta,s)$.\\ Since $(z,\theta)\mapsto g_i(z,\theta)$ is
continuous on $[-1,3]\times (-1/2,1/2)\!\setminus\!\{0\}$ we
finally deduce that for any $i \in \{1,2\}$, $\theta \mapsto
\int_{\theta s_1}^{\beta(\theta,s)} g_i(\rho,\theta,s)d\rho$ is
differentiable on $(-1/2,1/2)\!\setminus\!\{0\}$ and,
\begin{eqnarray*}
\sum_{i=1}^2\frac{\partial}{\partial \theta}\int_{\theta
s_1}^{\beta(\theta,s)} g_i(\rho,\theta,s) d\rho
&=&
\frac{\theta(s_1^2-1)}{\sqrt{1-\theta^2+\theta^2s_1^2}}\\
&&
\times\Big(-(n-3)s_1(\sqrt{1-\theta^2+\theta^2s_1^2}-\theta
s_1)^{n-4}\\
&&
\quad
+\,\alpha\theta(1-s_1^2)(\sqrt{1-\theta^2+\theta^2s_1^2}-\theta
s_1\big)^{n-3}\Big)\\
&+&
\sum_{i=1}^2\int_{\theta s_1}^{\beta(\theta,s)} \frac{\partial^2
g_i}{\partial^2 \theta}(\rho,\theta,s)d\rho.
\end{eqnarray*}
We deduce that $\theta \mapsto \frac{\partial F}{\partial \theta}$
is differentiable in $(-1/2,1/2)\!\setminus\!\{0\}$. Moreover, we
see that the map,
\begin{eqnarray*}
\theta \mapsto \lambda (\theta,s)&=&
\frac{\theta(s_1^2-1)}{\sqrt{1-\theta^2+\theta^2s_1^2}}\Big(
-(n-3)s_1(\sqrt{1-\theta^2+\theta^2s_1^2}-\theta s_1)^{n-4}\\
&&
\hspace{3cm}
+\,\alpha\theta(1-s_1^2)(\sqrt{1-\theta^2+\theta^2s_1^2}-\theta
s_1\big)^{n-3}\Big)
\end{eqnarray*}
is indefinitely differentiable on $(-1/2,1/2)\times
\mathbb{S}^{n-1}$. Then,
by $(1.1),(1.2),(1.8),(1.7),(1.3)$ and $(A)$, for any $a,b \in
(0,1/2)$, $a<b$ there exists constants $K_{1,n,ab,\alpha},
K_{2,n,ab,\alpha}, K_{3,n,ab,\alpha}$ so that, for any $|\theta|
\in (a,b)$, for any $s \in \mathbb{S}^{n-1}\setminus \{(\pm
1,0,\cdots,0)\}$,
\[
\left|\frac{\partial^2 HF}{\partial \theta^2}(\theta,s)\right|
\leq
K_{1,n,ab,\alpha}(1-s_1^2)^{\frac{\alpha}{2}}+K_{2,n,ab,\alpha}(1-s_1^2)^{\frac{3-n}{2}}+
K_{3,n,ab,\alpha}.
\]
We deduce by Lebesgue Theorem that the map $\theta \mapsto
E(\theta)$ is two times differentiable on
$(-1/2,1/2)\!\setminus\!\{0\}$ and,
\[
\frac{d^2 G}{d \theta^2}(\theta)=
\int_{\mathbb{S}^{n-1}}\frac{\partial^2 HF}{\partial
\theta^2}(\theta,s))ds.
\]
\end{dem}
\begin{lem}
If $5-n> \alpha > 4-n$, the map $\theta \mapsto G(\theta)$ is two
times differentiable at $0$.
\end{lem}
\begin{dem}
Suppose that $\alpha \in (4-n,5-n)$. As in  Lemma 1.5, we can see
that,
\[
\int_{\mathbb{S}^{n-1}}\frac{\partial^2 H}{\partial \theta^2
}(\theta,s)F(\theta,s)ds\mathop{\longrightarrow}_{\theta
\rightarrow 0}
\int_{\mathbb{S}^{n-1}}\frac{1}{2}\frac{(2n-3)s_1^2-(n-2)}{n-2+\alpha}
ds=\frac{-n^2+4n-3}{2n(n-2+\alpha)}|\mathbb{S}^{n-1}|,
\]
\[
\int_{\mathbb{S}^{n-1}}\frac{\partial H}{\partial \theta
}(\theta,s)\eta(\theta,s)ds \mathop{\longrightarrow}_{\theta
\rightarrow 0} 0,
\]
\[
\int_{\mathbb{S}^{n-1}}\frac{\partial H}{\partial \theta
}(\theta,s)\int_{\theta
s_1}^{\beta(\theta,s)}g_1(\rho,\theta,s)d\rho
\mathop{\longrightarrow}_{\theta \rightarrow
0}\int_{\mathbb{S}^{n-1}}\frac{2(n-3)(n-1)}{n-3+\alpha}s_1^2=
\frac{2(n-3)(n-1)}{n(n-3+\alpha)}|\mathbb{S}^{n-1}|,
\]
\[
\int_{\mathbb{S}^{n-1}}\frac{\partial H}{\partial \theta
}(\theta,s)\int_{\theta
s_1}^{\beta(\theta,s)}g_2(\rho,\theta,s)d\rho
\mathop{\longrightarrow}_{\theta \rightarrow 0}0,
\]
\[
\int_{\mathbb{S}^{n-1}}H(\theta,s)\frac{\partial \eta}{\partial
\theta }(\theta,s)ds \mathop{\longrightarrow}_{\theta \rightarrow
0} \int_{\mathbb{S}^{n-1}}(n-1)(s_1^2-1)ds=
\frac{-(n-1)^2}{n}|\mathbb{S}^{n-1}|,
\]
and
\[
\int_{\mathbb{S}^{n-1}}H(\theta,s)\lambda(\theta,s)ds
\mathop{\longrightarrow}_{\theta \rightarrow 0} 0.
\]
As in Lemma 1.5, we set $\rho=\sqrt{1-s_1^2}\theta t$ if
$\theta>0$. Hence,
\begin{eqnarray*}
\int_{\theta s_1}^{\beta(\theta, s)}g_{11}(\rho,\theta,s)d\rho&=&
(n\!-\!3)(n\!-\!4)s_1^2(1-s_1^2)^{\frac{1+\alpha}{2}}
\theta^{n-4+\alpha}J(n-5,\frac{\alpha}{2}).
\end{eqnarray*}
\begin{eqnarray*}
\int_{\theta s_1}^{\beta(\theta s)}g_{12}(\rho,\theta,s)d\rho
&=&
-2\alpha(n\!-\!3)s_1(1-s_1^2)^{\frac{1+\alpha}{2}}\theta^{n-4+\alpha}J(n-4,\frac{\alpha}{2}-1)
\end{eqnarray*}
\begin{eqnarray*}
\int_{\theta s_1}^{\beta(\theta s)}g_{21}(\rho,\theta,s)d\rho
&=&
\alpha(\alpha-2)(1-s_1^2)^{\frac{1+\alpha}{2}}\theta^{n-4+\alpha}J(n-3,\frac{\alpha}{2}-2)
\end{eqnarray*}
\begin{eqnarray*}
\int_{\theta s_1}^{\beta(\theta s)}g_{22}(\rho,\theta,s)d\rho
&=&
\alpha
(1-s_1^2)^{\frac{1+\alpha}{2}}\theta^{n-4+\alpha}J(n-3,\frac{\alpha}{2}-1)
\end{eqnarray*}
Since $\alpha \in (4-n,5-n)$, the integrals
$J(n-5,\frac{\alpha}{2})$ and $J(n-3,\frac{\alpha}{2}-1)$
are infinite and we have,
\[
J(n-5,\frac{\alpha}{2}) \mathop{\sim}_{0}
\frac{(1\!-\!s_1^2)^{\frac{-1-\alpha}{2}}\theta^{-n-\alpha+4}}{n-4+\alpha},\,
J(n-3,\frac{\alpha}{2}-1) \mathop{\sim}_{0}
\frac{(1\!-\!s_1^2)^{\frac{1-\alpha}{2}}\theta^{-n-\alpha+4}}{n-4+\alpha}
\]
And the integrals
$J(n-4,\frac{\alpha}{2}-1)$ and $J(n-3,\frac{\alpha}{2}-2)$
are finite. Then,
\[
\int_{\theta s_1}^{\beta(\theta, s)}g_{11}(\rho,\theta,s)d\rho
\mathop{\longrightarrow}_{\theta \mapsto 0}
\frac{(n-3)(n-4)s_1^2}{n-4+\alpha},\,
\int_{\theta s_1}^{\beta(\theta, s)}g_{22}(\rho,\theta,s)d\rho
\mathop{\longrightarrow}_{\theta \mapsto 0}
\frac{\alpha(1\!-\!s_1^2)}{n-4+\alpha}.
\]
\[
\int_{\theta s_1}^{\beta(\theta, s)}g_{12}(\rho,\theta,s)d\rho
\mathop{\longrightarrow}_{\theta \mapsto 0} 0,\,
\int_{\theta s_1}^{\beta(\theta, s)}g_{21}(\rho,\theta,s)d\rho
\mathop{\longrightarrow}_{\theta \mapsto 0} 0
\]
Moreover, we can see that, for any $i,j \in \{1,2\}$, for any $
(\theta,s) \in (-1/2,1/2) \times
\mathbb{S}^{n-1}\!\setminus\!\{(\pm 1,0,\cdots,0)\}$,
\begin{eqnarray*}
H(\theta,s)\int_{\theta s_1}^{\beta(\theta
s)}g_{ij}(\rho,\theta,s)d\rho
&\leq&
C_{n,\alpha}(1-s_1^2)^{\frac{5-n}{2}} +
D_{n,\alpha}(1-s_1^2)^{\frac{\alpha+1}{2}}.
\end{eqnarray*}
where $C_{n,\alpha}$ and $D_{n,\alpha}$ are two constants
independent of $\theta$. By Lebesgue Theorem we deduce that,
\begin{eqnarray*}
\int_{\mathbb{S}^{n-1}}H(\theta,s)\frac{\partial^2 F}{\partial
\theta^2 }(\theta,s)ds \mathop{\longrightarrow}_{\theta
\rightarrow 0}
&&
\frac{-(n-1)^2}{n}|\mathbb{S}^{n-1}|\\
&&
\quad +\,(n-1) \frac{(n-3)(n-4)+\alpha
(n-1)}{n(n-4+\alpha)}|\mathbb{S}^{n-1}|.
\end{eqnarray*}
By Lemmas $1.1, 1.2, 1.3$, $\theta \mapsto G(\theta) \in
\mathcal{C}^1((-1/2,1/2),\mathbb{R})$ and is two times
differentiable on $(-1/2,1/2)\!\setminus\!\{0\}$. Furthermore,
when $\alpha \in (4-n,5-n)$, as  the limit of $\frac{d^2 G}{d
\theta^2}(\theta)$ exists as $\theta \rightarrow 0$ , we have $
\theta \mapsto G(\theta)$ is two times differentiable on
$(-1/2,1/2)$.
\end{dem}
\newline
\textbf{Proof of ii).}
Assume that $\alpha \in (4-n,5-n)$, by Lemma 1.1, 1.2, 1.3, 1.4,
we have,
\[
G(\theta)= G(0)+ \frac{1}{2}\frac{d^2 G}{d \theta^2}(0)+
o(\theta^2).
\]
Furthermore we have,
\begin{eqnarray*}
\frac{d^2 G}{d \theta^2}(0)
&=&
\frac{-n^2+4n-3}{2n(n-2+\alpha)}|\mathbb{S}^{n-1}|
+
\frac{2(n-3)(n-1)}{n(n-3+\alpha)}|\mathbb{S}^{n-1}|\\
&+&
\frac{-(n-1)^2}{n}|\mathbb{S}^{n-1}|
\quad +\, (n-1)\frac{(n-3)(n-4)+\alpha
(n-1)}{n(n-4+\alpha)}|\mathbb{S}^{n-1}|.
\end{eqnarray*}
We have, for any $n \geq 6$.
\[
(n-3)(n-4)+\alpha (n-1) \mathop{\longrightarrow}_{\alpha \mapsto
4-n} -2(n-4) <0.
\]
Then,
\[
\frac{(n-3)(n-4)+\alpha
(n-1)}{n(n-4+\alpha)}\mathop{\longrightarrow}_{\alpha \mapsto_{>}
4-n}-\infty,\, \text{and}\, \frac{d^2 G}{d^2
\theta}(0)\mathop{\longrightarrow}_{\alpha \mapsto_{>}
4-n}-\infty.
\]
Hence, there is $\alpha_0$ such that, for any $\alpha \in
(4-n,\alpha_0)$, $G(\theta) < G(0)$ for $\theta$ sufficiently
small, that is,
\[
G(\theta)=E_{2,r^{\alpha}}(u_a)=\int_ {\mathbf{B}^n} r^{\alpha}
\|\nabla u_a\|^2 dx < G(0)=\int_ {\mathbf{B}^n} r^{\alpha}
\|\nabla u_0\|^2 dx. \quad \quad \quad \blacksquare
\]

\begin{thebibliography}{}
{\small\bibitem{AL}\textsf{M. Avellaneda, F.-H. Lin},\,
Null-Lagrangians and minimizing $\int |\nabla u|^p$ , \textit{C.
R. Acad. Sci. Paris}, \textbf{306}(1988), 355-358.
%
%
\bibitem{BCL}\textsf{H. Brezis, J.-M. Coron, E.-H. Lieb}, Harmonic Maps with
Defects, \textit{Commun. Math. Phys.}, \textbf{107} (1986),
649-705.
%
%
\bibitem{C}\textsf{B. Chen},\, Singularities of $p$-harmonic mappings, \textit{Thesis, University of Minnesota},
(1989).
%
%
\bibitem{CHe}\textsf{J.-M. Coron, F. Helein}, Harmonic diffeomorphisms, minimizing
harmonic maps and rotational symmetry, \textit{Compositio Math.},
\textbf{69} (1989), 175-228.
%
%
\bibitem{CG}\textsf{J.-M. Coron, R. Gulliver}, Minimizing $p$-harmonic maps into spheres, \textit{J. reine angew. Math.},
\textbf{401} (1989), 82-100.
%
%
\bibitem{EJ}\textsf{A. El Soufi, A. Jeune}, Indice de Morse des applications $p$-harmoniques,
\textit{Ann. Inst. Henri Poincaré}, \textbf{13}(1996), 229-250.
%
%
\bibitem{ES}\textsf{A. El Soufi, E. Sandier},\, $p$-harmonic
diffeomorphisms, \textit{Calc. of Var.}, \textbf{6}(1998),
161-169.
%
%
\bibitem{CH}\textsf{B. Chen, R. Hardt}, Prescribing singularities for $p$-harmonic mappings,
\textit{Indiana University Math. J.}, \textbf{44}(1995), 575-601.
%
%
bibitem{HL1}\textsf{R. Hardt, F.-H. Lin}, Mapping minimizing the
$L^p$ norm of the gradient, \textit{Comm. P.A.M.},
\textbf{15}(1987), 555-588.
%
%
\bibitem{
HL2}\textsf{R. Hardt, F.-H. Lin}, Singularities for $p$-energy
minimizing unit vectorfields on planar domains , \textit{Calculus
of Variations and Partial Differential Equations},
\textbf{3}(1995), 311-341.
%
%
\bibitem{HLW1}
\textsf{R. Hardt, F.-H. Lin, C.-Y. Wang}, The $p$-energy
minimality of $\frac{x}{\|x\|}$, \textit{Communications in
analysis and geometry}, \textbf{6}(1998), 141-152.
%
%
\bibitem{
HLW2}\textsf{R. Hardt, F.-H. Lin, C.-Y. Wang}, Singularities of
$p$-Energy Minimizing Maps, , \textit{Comm. P.A.M.},
\textbf{50}(1997), 399-447.
%
%
\bibitem{He}\textsf{F. Helein}, Harmonic diffeomorphisms between an open subset of $ \mathbb{R}^3$
and a Riemannian manifold, \textit{C. R. Acad. Sci. Paris},
\textbf{308}(1989), 237-240.
%
%
\bibitem{HKW}
\textsf{S. Hildebrandt, H. Kaul, K.-O. Wildman}, An existence
theorem for harmonic mappings of Riemannian manifold, \textit{Acta
Math.}, \textbf{138}(1977), 1-16.
%
%
\bibitem{Ho1}\textsf{M.-C. Hong,}, On the Jäger-Kaul theorem concerning harmonic maps,
\textit{Ann. Inst. Poincaré, Analyse non-linéaire}, \textbf{17}
(2000), 35-46.
%
%
\bibitem{Ho2}\textsf{M.-C. Hong}, On the minimality of the $p$-harmonic map
$\frac{x}{\|x\|}:\mathbf{B}^n \rightarrow \mathbf{S}^{n-1}$,
\textit{Calc. Var.}, \textbf{13} (2001), 459-468.
%
%
\bibitem{JK}\textsf{W. Jager, H. Kaul}, Rotationally symmetric harmonic maps from a ball
into a sphere and the regularity problem for weak solutions of
elliptic systems, \textit{J.Reine Angew. Math.},
\textbf{343}(1983), 146-161.
%
%
\bibitem{L}\textsf{F.-H. Lin}, Une remarque sur l'application $x/\|x\|$, \textit{C.R. Acad. Sci.
Paris} \textbf{305}(1987), 529-531.
%
%
\bibitem{SU1}\textsf{R. Shoen, K. Uhlenbeck}, A regularity theory for harmonic maps \textit{J. Differential Geom.}.
\textbf{12}(1982), 307-335.
%
%
\bibitem{SU2}\textsf{R. Shoen, K. Uhlenbeck}, Boundary theory and the Dirichlet
problem for harmonic maps, \textit{J. Differential Geom.},
\textbf{18}(1983), 253-268.}
%
%
%
%
\bibitem{W} \textsf{C.Wang}, Minimality and perturbation
of singularities for certain $p$-harmonic maps., \textit{Indania
Univ. Math.J.}, \textbf{47}(1998), 725-740.
%
%
\
\end{thebibliography}
\end{document}